\documentclass[12pt,oneside,reqno]{amsart}

\usepackage[T1]{fontenc}
\usepackage{lmodern}

\usepackage{amssymb,amsfonts,amsmath,amsthm,mathtools,cancel}

\usepackage{graphicx,float,wrapfig}

\usepackage{caption,subcaption}

\usepackage{geometry}

\usepackage{xcolor}

\usepackage{tabularray}

\usepackage{tikz}
\usetikzlibrary{cd,graphs}
\usepackage[all]{xy}

\usepackage{enumitem}

\usepackage{chngcntr}

\usepackage{url}

\usepackage[hidelinks, bookmarks=false]{hyperref}
\newcommand{\DOI}[2]{{\fontfamily{lmtt}\bfseries\href{#1}{#2}}}

\def\qed{\ifhmode\textqed\fi
	\ifmmode\ifinner\hfill\quad\qedsymbol\else\dispqed\fi\fi}
\def\textqed{\unskip\nobreak\penalty50
	\hskip2em\hbox{}\nobreak\hfill\qedsymbol
	\parfillskip=0pt \finalhyphendemerits=0}
\def\dispqed{\rlap{\qquad\qedsymbol}}

\setlist[enumerate]{%
	label={\normalfont\arabic*.},
	left={\parindent},
	itemsep={5pt}%
}%
\newlist{alphanumerate}{enumerate}{10}
\setlist[alphanumerate]{%
	label={\normalfont\alph*)},%
	left={\parindent},%
	itemsep={5pt}%
}%

\newcommand{\MBB}{\mathbb}
	\def\ZZ{\MBB Z} 

\let\Sect=\bigcap

\let\epsilon=\varepsilon
\let\phi=\varphi
\let\kappa=\varkappa

\def\Implies{\ifmmode\Longrightarrow \else
	\unskip${}\Longrightarrow{}$\ignorespaces\fi}
\def\implies{\ifmmode\Rightarrow \else
	\unskip${}\Rightarrow{}$\ignorespaces\fi}
\def\iff{\ifmmode\Longleftrightarrow \else
	\unskip${}\Longleftrightarrow{}$\ignorespaces\fi}

\let\:=\colon

\DeclareMathOperator{\depth}{depth}

\DeclareMathOperator{\Coker}{Coker}

\DeclareMathOperator{\Ext}{Ext}

\theoremstyle{plain}
	\newtheorem{Theorem}{Theorem}[section]
	
	\newtheorem{Proposition}[Theorem]{Proposition}

\theoremstyle{definition}
	\newtheorem{Definition}[Theorem]{Definition}
	\newtheorem{Example}[Theorem]{Example}
	
	\newtheorem{Remark}[Theorem]{Remark}
\def\CM{{Cohen-Macaulay}}
\def\SCM{{sequentially Cohen-Macaulay}}
\def\CCM{{canonically Cohen-Macaulay}}
\def\FCM{{Cohen-Macaulay filtered}}

\begin{document}
\title{SCMAlgebras: a Macaulay2 package for sequentially Cohen-Macaulayness}
\author{Ernesto Lax}

\address{Ernesto Lax, Department of mathematics and computer sciences, physics and earth sciences, University of Messina, Viale Ferdinando Stagno d'Alcontres 31, 98166 Messina, Italy}
\email{erlax@unime.it}

\thanks{{\tt SCMAlgebras v1.2}}

\begin{abstract}
  We introduce the \textit{Macaulay2} package \texttt{SCMAlgebras}. It provides functions for computing the modules of deficiency and the filter ideals, in order to check whether a module or an ideal is {\SCM}. The package also implements routines for studying other {\CM} type conditions and unmixedness. After recalling the basic algebraic notions and results, the main features of the package are described through examples.
\end{abstract}

\maketitle

\section{Introduction}\label{sec:intro}
The notion of {\SCM} modules was originally introduced by Stanley \cite{Stan_book} in the graded setting, as an algebraic counterpart to the concept of nonpure shellable complexes, which had been developed by Björner and Wachs \cite{BW} in the context of combinatorics and discrete geometry. Later on, Schenzel \cite{S} independently introduced the notion of {\FCM} modules and proved the equivalence of this notion to the {\SCM} property. Since then, the theory of {\SCM} modules has garnered significant attention from many authors \cite{ABG,AM,CD,CDSS,D,ERT2022b,G,HTYZN,HS,LRR,SZ,VTV,W}, particularly due to its deep connections with the study of algebraic properties of combinatorial objects such as simplicial complexes and graphs.

In this paper, we present a new {\em Macaulay2} \cite{GDS} package, {\tt SCMAlgebras} \cite{L_pack}, designed to check whether a given module or ideal satisfies the {\SCM} property. This package implements key theoretical results, such as Theorem \ref{thm:SCM_DefModChar}, which characterizes {\SCM} modules in terms of their modules of deficiency, and Proposition \ref{prop:SCM-char-goodarzi} due to Goodarzi \cite{G}, which provides an efficient criterion for verifying the {\SCM} property for homogeneous ideal. To support these results, the package introduces computational tools for calculating the modules of deficiency (see Definition \ref{def:ModulDef}) and generating filter ideals (see Definition \ref{def:FilterIdeal}), both of which are fundamental for testing the {\SCM} property of the input objects. Beyond this primary goal, these constructions also enable the computation of related invariants and the study of further properties. In particular, the filter ideals provide an effective approach for checking the unmixedness of an ideal (Proposition \ref{prop:unmixedness_char}) \emph{via} the notion of unmixed layer, while the package also includes routines for investigating {\CM} type conditions. In this way, {\tt SCMAlgebras} provides a broader computational framework for the study of modules and ideals, extending its applicability beyond the verification of the {\SCM} property.

These tools extend the functionality of the software, playing as a starting point for future development in the manipulation of {\SCM} property in both theoretical and applied contexts.

\section{Mathematical background}\label{sec:1:SCM_theory}
In this section we recall the notion of {\SCM} modules and the main characterizations, that are implemented in the package. In the following we will work in the graded case. Let $S=K[x_1,\ldots,x_n]$ be the standard graded polynomial ring in $n$ variables over a field $K$.\par

\begin{Definition}\label{def:SCM_definition}
    A finitely generated graded $S$-module $M$ is called a {\em\SCM} module if there exist a finite filtration of submodules of $M$
    \begin{equation*}
        0 = M_0 \subset M_1 \subset \cdots \subset M_{d-1} \subset M_r = M   
    \end{equation*}
    such that
    \begin{alphanumerate}
        \item each quotient $M_i/M_{i-1}$ is {\CM},
        \item $\dim M_i/M_{i-1} < \dim M_{i+1}/M_i$ for all $i$.
    \end{alphanumerate}
\end{Definition}

A filtration satisfying the above conditions is called a {\em{\SCM} filtration} of $M$. It is easy to see that if the {\SCM} filtration exists, then it is unique.\par

\begin{Remark}\label{rmk:CM->SCM}
    Any {\CM} module is {\SCM}, by considering the {\SCM} filtration $0=M_0 \subset M_1 = M$. Hence, one has the implication\par
    \noindent
    \begin{center}
        \begin{tblr}{colspec={ccc},rows={m},colsep=.4cm,width=\linewidth,rowsep=.5cm}
    	{\CM} & {$\Implies$} & {sequentially\\{\CM}}
        \end{tblr}
    \end{center}
    but the converse is not true in general.
\end{Remark}


Proving the {\SCM} property through the construction of a {\SCM} filtration is not always an easy task. The definition itself is not operative and it is not very efficient when trying to implement it on algorithms. In order to demonstrate this condition effectively and one must rely on algebraic or combinatorial characterization.

In the ladder, two of the main characterizations of the sequentially Cohen-Macaulay property are recalled. The first one involves the notion of the modules of deficiency. 

\begin{Definition}\label{def:ModulDef}
  Let $M$ be a finitely generated graded $S$-module and let $d=\dim M$. For any $i\in\ZZ$ the {\em $i$th module of deficiency} of $M$ is 
  \[
  \omega^{i}(M)=\Ext^{n-i}_S(M,S(-n)).
  \]
  If $i=d$, the module $\omega^{d}(M)=\omega(M)$ is the {\em canonical module} of $M$.
\end{Definition}

\begin{Remark}\label{rem:Grot-Matlis}
    The modules of deficiency were introduced and studied by Schenzel in \cite{S_book1} in the local case. By duality (in the graded setting), one has that $\omega^i(M)$ is the Matlis dual of the $i$th local cohomology module $H^i(M)$ of $M$. Let $t=\depth M$.
    
    From a well--known result of Grothendieck (see \cite{BH}), it holds that $H^i(M)=0$ for $i<t$ or $i>d$. Moreover, $H^i(M)\neq 0$ for $i=t,d$. It follows that $M$ is {\CM} if and only if $H^i(M)=0$ for all $i\neq d$. Hence, one has $\omega^i(M) = 0$ for $i<t$ or $i>d$, $\omega^t(M)\neq 0$ and $\omega(M)=\omega(M)\neq 0$ and $M$ is a {\CM} module if and only if $\omega^i(M) = 0$ for all $i\neq d$. In such case, the canonical module $\omega^d(M)=\omega(M)$ is also {\CM}. Hence, the modules of deficiency of $M$ measure how far is $M$ from being a {\CM} module.
\end{Remark}

\begin{Definition}
    A module $M$ with {\CM} canonical module $\omega(M)$ is called a {\em\CCM} module.
\end{Definition}

The next classical result characterizes {\SCM} modules {\em via} the modules of deficiency (see \cite{CDSS}).

\begin{Theorem}\label{thm:SCM_DefModChar}
    Let $M$ be a finitely generated $S$-module, with $d=\dim M$. Then, the followings are equivalent:
    \begin{enumerate}
	\item $M$ is {\SCM};\label{thm:SCM_DefModChar:(1)}
        \item for all $0\leq i \leq d$ the $i$th module of deficiency $\omega^i(M)$ is either zero or an $i$-dimensional {\CM} module;\label{thm:SCM_DefModChar:(2)}
        \item for all $0\leq i < d$ the $i$th module of deficiency $\omega^i(M)$ is either zero or an $i$-dimensional {\CM} module. \label{thm:SCM_DefModChar:(3)}
    \end{enumerate}
\end{Theorem}

The equivalence \ref{thm:SCM_DefModChar:(1)}\ $\iff$\ \ref{thm:SCM_DefModChar:(2)} in the above theorem was announced without proof in \cite[Theorem 2.11]{Stan_book}, but the author mentioned a spectral sequence argument due to Peskine.

\begin{Remark}\label{rmk:SCM->CCM}
    If $M$ is a {\SCM} module, then by condition \ref{thm:SCM_DefModChar:(2)} of Theorem \ref{thm:SCM_DefModChar} $M$ has a {\CM} canonical module. Hence, it follows that a {\SCM} module is {\CCM}, but the converse is not true in general.
\end{Remark}\par

Thus, merging the informations from Remark \ref{rmk:CM->SCM} and Remark \ref{rmk:SCM->CCM}, one has the following chain of implications:\par\medskip

\noindent
\begin{center}
    \begin{tblr}{colspec={ccccc},rows={m},colsep=.4cm,width=\linewidth,rowsep=.5cm}
	{\CM} & {$\Implies$} & {sequentially\\{\CM}} & {$\Implies$} & {canonically\\{\CM}}
    \end{tblr}
\end{center}\medskip
with the reverse implications not true in general.\par\bigskip

Although Theorem \ref{thm:SCM_DefModChar} remains useful in many situations (see \cite[Theorem 1.5]{LRR},\cite{SZ}), it also presents two significant challenges. Firstly, to determine which modules of deficiency vanish is often a difficult task (see \cite{W}). Secondly, proving the {\CM}ness of the non--vanishing modules is typically non--trivial. These difficulties underscore the complexity of applying the theorem in practice, particularly when dealing with more complicated algebraic structures.\par\bigskip\bigskip

In case we work with a homogeneous ideal $I\subset S$ it is possible to characterize the {\SCM} property by using new ideals, namely the filter ideals, which are constructed from the primary decomposition of $I$.

\begin{Definition}\label{def:FilterIdeal}
	Let $I\subset S$ be a homogeneous ideal, with $d=\dim S/I$, and let $I=\Sect_{j=1}^r Q_j$ be the minimal primary decomposition of $I$. For all $1\leq j\leq r$, let $P_j = \sqrt{Q_j}$ be the radical of $Q_j$. For all $-1\leq i\leq d$, the {\em $i$th filter ideal} of $I$ is
	\[
	I^{<i>} = \Sect_{\dim S/{P_j}>i} Q_{j},
	\]
	where $I^{<-1>}=I$ and $I^{<d>}=S$. 
\end{Definition}

Moreover, we define the {\em minimum dimension} of $I$ as the minimum integer such that $I^{<i>}\neq I$. By definition, it is clear that minimum dimension is at least $-1$.

For a homogeneous ideal $I\subset S$, one has the filtration
\[
I = I^{<-1>} \subseteq I^{<0>} \subseteq I^{<1>} \subseteq \ldots \subseteq I^{<d-1>} \subseteq I^{<d>} = S.
\]
Moreover, rewriting the above sequence of inclusions modulo $I$, we obtain
\[
0 \subseteq I^{<0>}/I \subseteq I^{<1>}/I \subseteq \ldots \subseteq I^{<d-1>}/I \subseteq S/I,
\]
which is called the {\em dimension filtration} of $S/I$. One can show (see \cite[Lemma 1]{G}) that such filtration does not depend on the primary components, therefore it is unique. Note that, the filter ideals $I^{<i>}$ are also independent from the primary decomposition.\par

By the above construction, one has the following characterization of the {\SCM} property for $S/I$ due to Goodarzi (see \cite[Proposition 16]{G}).

\begin{Proposition}\label{prop:SCM-char-goodarzi}
    Let $I\subset S$ be a homogeneous ideal and suppose $d=\dim S/I$. Then, $S/I$ is {\SCM} if and only if
    \begin{equation}\label{eq:SCM-char-goodarzi}
    \depth S/{I^{<i>}} \geq i+1,\qquad\text{for all $0\leq i < d$}.
    \end{equation}
\end{Proposition}

Goodarzi's characterization provides an efficient criterion to verify the {\SCM} property for a homogeneous ideal. It has been widely used, for instance, to study the sequential {\CM}ness of the binomial edge ideals of various classes of graphs, such as closed graphs (see \cite{ERT2022b}), block graphs, cycles, and cone graphs (see \cite{LRR}), playing a crucial role in the proofs. In particular, for binomial edge ideals, this result also permitted in \cite{LRR} the investigation of situations that were difficult to describe using Theorem \ref{thm:SCM_DefModChar}. This breakthrough allowed for a deeper understanding of cases previously out of reach, thus highlighting the power of Goodarzi's approach.\par\bigskip

We end this section with another useful application of filter ideals.

\begin{Definition}\label{def:UnmixedLayer}
	Let $I\subset S$ be a homogeneous ideal, with $d=\dim S/I$. For all $0\leq i\leq d$, the {\em $i$th unmixed layer} of $S/I$ is the $S$-module $I^{<i>}/I^{<i-1>}$ and it is denoted by $U_i(S/I)$.
\end{Definition}

The name unmixed layer is motivated by \cite[Corollary 2.3]{S}, in fact each module $U_i(S/I)$ is unmixed, that is for all $i$ the associated primes of $U_i(S/I)$ have the same height. The next result (see \cite[Lemma 3]{G}) gives a characterization of the unmixedness of the ideal $I$.

\begin{Proposition}\label{prop:unmixedness_char}
    Let $I\subset S$ be a homogeneous ideal, with $d=\dim S/I$. Then $S/I$ is unmixed if and only if $U_i(S/I)=0$ for all $1 \leq i < d$ and $U_d(S/I)=S/I$.
\end{Proposition}\pagebreak

\section{The package}\label{sec:Package}
In this section we describe the functions of the package {\tt SCMAlgebras} \cite{L_pack} by showing some examples of usage and algorithms.

\begin{Example}\label{ex:module_matrix}
    Let $S = K[x_1,x_2,x_3,x_4,x_5]$ and let $M$ be the $S$-module defined as the $\Coker$ of the map represented by the matrix 
    $\left(\begin{smallmatrix}
    	x_1x_2 & x_3x_4 & 0 & 0\\
    	0 & x_1x_5 & x_2x_4 &0
    \end{smallmatrix}\right)$. 
    Such a module has $\depth M = 3$ and $\dim M = 4$. We use the {\tt deficiencyModule} function to compute the modules of deficiency of $M$.\par
    
    \verb/  i1 :  loadPackage "SCMAlgebras";/
    
    \verb/  i2 :  S=QQ[x_1..x_5];/

    \verb/  i3 :  M=coker matrix{{x_1*x_2,x_3*x_4,0,0},{0,x_1*x_5,x_2*x_4,0}}/
    
    \verb/  o3 =  cokernel | x x   x x   0     0 | /\par\vspace{-5pt}
    \verb/                    1 2   3 4             /\par
    \verb/                 | 0     x x   x x   0 | /\par\vspace{-5pt}
    \verb/                          1 5   2 4       /  

    \verb/                               2 /\par\vspace{-5pt}
    \verb/  o3 :  S-module, quotient of S  /\medskip    
    
    \verb/  i4 :  deficiencyModule(M,1) /
    
    \verb/  o4 =  0 /

    \verb/  o4 :  S-module /
    
    \verb/  i5 :  deficiencyModule(M,3) /
    
    \verb/                    2       2           /\par\vspace{-5pt}
    \verb/  o5 =  cokernel | x x   x x   x x x  | /\par\vspace{-5pt}
    \verb/                    1 5   3 4   1 2 4   /

    \verb/                               1 /\par\vspace{-5pt}
    \verb/  o5 :  S-module, quotient of S  /\bigskip

    The function {\tt deficiencyModule(M,i)} uses Definition \ref{def:ModulDef} and computes the module $\Ext^{n-i}_S(M,S(-n))$. According to Remark \ref{rem:Grot-Matlis}, in our case, the function returns the zero module whenever $i<3$ or $i>4$, and the only non vanishing modules of deficiency are $\omega^3(M)$ and $\omega^4(M)=\omega(M)$, {\em i.e.} the canonical module of $M$. Moreover, the function {\tt canonicalModule(M)} allows to compute the canonical module $\omega(M)$ directly, as the module of deficiency corresponding to dimension.\par\bigskip

    \verb/  i6 :  canonicalModule(M) == deficiencyModule(M,4) /

    \verb/  o6 =  true /\bigskip

    Hence, one can check if $M$ is {\SCM} using the {\tt isSCM} function, involving the modules of deficiency, as seen in Theorem \ref{thm:SCM_DefModChar}.\bigskip

    \verb/  i7 :  isSCM M /
    
    \verb/  o7 =  false /\bigskip
    
    Moreover, using the {\tt isCCM} function, one can check that $M$ is {\CCM}, by checking the Cohen-Macaulayness of its canonical module.\bigskip
    
    \verb/  i8 :  isCCM M /
        
    \verb/  o8 =  true /
\end{Example}

\begin{Remark}
	Note that, while in Example \ref{ex:module_matrix}, the functions {\tt deficiencyModule} and {\tt canonicalModule} are applied to a module, they can also be applied to an ideal.
\end{Remark}

\begin{Example}\label{ex:BinomialEdge}
  Consider the graph $G$ depicted below.
  
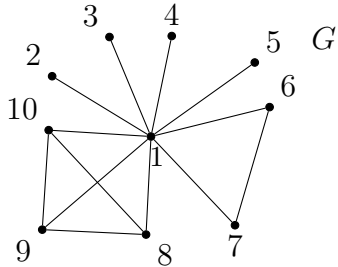
\begin{figure}[H]
    \centering
    \begin{tikzpicture}[scale=.65,vertices/.style={draw, fill=black, circle, inner sep=1pt}]
      \node[vertices] (v) at (0,0) {}; 
      \node[vertices] (1) at (-2, 1.23) {}; 
      \node[vertices] (2) at (-0.84, 2.03) {}; 
      \node[vertices] (3) at (0.42, 2.05) {};	
      \node[vertices] (4) at (2.1, 1.51) {}; 
      \node[vertices] (5) at (2.4,0.6) {}; 
      \node[vertices] (6) at (1.7, -1.81) {};		
      \node[vertices] (7) at (-0.1,-2) {}; 
      \node[vertices] (8) at (-2.2,-1.9) {};
      \node[vertices] (9) at (-2.07, 0.13) {};
      \node at (.12,-.41) {\small$1$};
      \node at (1) [above left] {$2$};
      \node at (2) [above left] {$3$};  
      \node at (3) [above] {$4$};  
      \node at (4) [above right] {$5$};  
      \node at (5) [above right] {$6$}; 
      \node at (6) [below] {$7$};
      \node at (7) [below right] {$8$};
      \node at (8) [below left] {$9$};
      \node at (9) [above left] {$10$};
    	
      \foreach \k in {1,...,9}{
          \draw (v)--(\k);
      };
    	
      \draw (5)--(6);
    	
      \draw (7)--(8);
      \draw (7)--(9);
      \draw (8)--(9);
    	
  	  \node at (3.5,2) {$G$};
    \end{tikzpicture}
  \end{figure}
    
  Consider the ideal $J_G$ in the polynomial ring $S=K[x_1,\ldots,x_{10},y_1,\ldots,y_{10}]$ generated by all the binomials $f_{ij} = x_iy_j - x_jy_i$, such that $\{i,j\}$ with $i<j$ is an edge of $G$, namely the {\em binomial edge ideal} (see \cite{HHHKR,O}) of $G$.\par\bigskip
  
  \verb/  i1 :  loadPackage "SCMAlgebras"; /   
  
  \verb/  i2 :  S = QQ[x_1..x_10,y_1..y_10]; /
  
  \verb/  i3 :  E = {{1,2},{1,3},{1,4},{1,5},{1,6},{1,7},{1,8},{1,9},{1,10},/
  
  \verb/        {6,7},{8,9},{8,10},{9,10}}; /
  
  \verb/  i4 :  J=ideal(for e in E list x_(e#0)*y_(e#1)-x_(e#1)*y_(e#0)); /\par\bigskip

  The function {\tt filterIdeals} computes the filter ideals of $J_G$, according to Definition \ref{def:FilterIdeal}. Moreover, the function {\tt minimumDimension} returns the smallest integer such that the filter ideals differ from $J_G$.\medskip
  
  \verb/  i5 :  filterIdeal(J,13)/
          
  \verb/  o5 =  ideal (y , x , x  y  - x y  , x  y  - x y  , x y  - x y , /\par\vspace{-5pt}  
  \verb/                1   1   10 9    9 10   10 8    8 10   9 8    8 9  /\par
  \verb/        --------------------------------------------------------- /\par
  \verb/        x y  - x y ) /\par\vspace{-5pt}  
  \verb/         7 6    6 7 /\pagebreak

  \verb/  o5 :  Ideal of S/
  
  \verb/  i6 :  minimumDimension J/
  
  \verb/  o6 =  11 /
  
  \verb/  i7 :  filterIdeal(J,7) == J/
  
  \verb/  o7 =  true /
  
  \verb/  i8 :  filterIdeal(J,11) == J/
  
  \verb/  o8 =  false /\medskip

  Using the function {\tt unmixedLayer}, one can compute the unmixed layers of $J_G$, according to Definition \ref{def:UnmixedLayer}.\par\medskip
  
  \verb/  i9  :  unmixedLayer(J,13)/
  
  \verb/  o9  =  subquotient (| y   x   x  y  - x y    x  y  - x y      /\par\vspace{-5pt}
  \verb/                         1   1   10 9    9 10   10 8    8 10/\par
      
  \verb/         x y  - x y   x y  - x y  |, | y   x   x  y  - x y     /\par\vspace{-5pt}
  \verb/          9 8    8 9   7 6    6 7       1   1   10 9    9 10  /\par
  
  \verb/         x  y  - x y    x y  - x y   x y  - x y  |)  /\par\vspace{-5pt}
  \verb/          10 8    8 10   9 8    8 9   7 6    6 7    /\par
  
  \verb/                                1 /\par\vspace{-5pt}
  \verb/  o9  :  S-module, quotient of S  /\medskip

  To check the unmixedness of $J_G$, the function {\tt isUnmixed} can be applied to the ideal {\tt J}. In this case Proposition \ref{prop:unmixedness_char} is used.\par\medskip
  
  \verb/  i10 :  isUnmixed J /

  \verb/  o10 =  false /\medskip    

  Finally, one can check that $S/J_G$ is {\SCM} (as proved in \cite[Theorem 1.5]{LRR}) using the {\tt isSCM} function. This time we apply the function to the ideal {\tt J}, in such case Proposition \ref{prop:SCM-char-goodarzi} is used.\par\medskip
  
  \verb|  i11 :  isSCM J|
  
  \verb|  o11 :  true |
\end{Example}

To optimize the computations, the new {\tt SCMAlgebras v1.2} introduces a major refactor of its internal routines along with new structures specifically designed to handle primary decompositions efficiently.

In previous versions, iterative calls to {\tt filterIdeal} function forced recomputing the minimal primary decomposition of the ideal each time. To bypass this computational bottleneck, a new data type called {\tt PrimaryDataList} has been introduced. This type stores the relevant information from the primary decomposition of an ideal $I$ as a list of all pairs $\{Q, \dim(S/\sqrt{Q})\}$, where $Q$ is a primary component of $I$. The user can generate this data via the new constructor function {\tt getPrimaryData(I)} and pass it as an optional argument to {\tt filterIdeal}.

In the present release, this strategy is internally adopted by the functions {\tt isSCM} and {\tt unmixedLayer}, which call {\tt filterIdeal} in a {\tt for} loop. By computing and caching the {\tt PrimaryDataList} upfront, the primary decomposition is performed only once at the beginning of the process rather than at each iteration, reducing execution times and memory leaks.

\begin{Remark}
  At the end of Example \ref{ex:module_matrix} we applied the {\tt isSCM} function to {\tt J} and Goodarzi's characterization (Proposition \ref{prop:SCM-char-goodarzi}) was used to check the {\SCM} property for $S/J_G$. Note that, to make the check involving Theorem \ref{thm:SCM_DefModChar}, one can apply the function to \verb|S^1/J| instead of {\tt J},\par\bigskip

	\verb|  i10 :  isSCM(S^1/J)|
	
	\verb|  o10 :  true |\par\bigskip
	
	\noindent
	getting, as expected, the same result.
\end{Remark}

We conclude by listing the main functionalities available in the package {\tt SCMAlgebras}.\par
\begin{table}[H]
    \centering
    \SetTblrInner{columns={l},rows={.75cm,m},hline{1,2,Z}={1.2pt},row{1}={lightgray!75}}
    \begin{tblr}{Q[1] | Q[2]}
        {\bf Function/Types} & {\bf Description}\\
        {\tt PrimaryDataList} & {New type that stores primary decomposition datas used for computations}\\
        {\tt getPrimaryData(I)} & {Computes a {\tt PrimaryDataList} of the ideal 
        $I$}\\
        {\tt deficiencyModule(M,i)} & {Computes the $i$th module of deficiency $\omega^i(M)$ of $M$}\\
        {\tt canonicalModule(M)} & {Computes the canonical module $\omega(M)$ of $M$}\\
        {\tt minimumDimension(I)} & Computes the minimum dimension of $I$\\		
        {\tt filterIdeal(I,i)} & Computes the $i$th filter ideal $I^{<i>}$ of $I$\\
        {\tt unmixedLayer(I,i)} & Computes the $i$th unmixed layer $U_i(I)$ of $I$\\
        {\tt isUnmixed(I)} & Checks if $I$ is unmixed \\
        {\tt isSCM(M)} & Checks if $M$ is {\SCM}\\
        {\tt isCCM(M)} & Checks if $M$ is {\CCM}
    \end{tblr}
    \caption{List of functions of {\tt SCMAlgebras v1.2}.}
\end{table}\par

\section{Conclusion and perspective}\label{sec:Conclusion}
The examples and the algorithms presented in this work are part of version {\tt 1.2} of the {\em Macaulay2} \cite{GDS} package {\tt SCMAlgebras} \cite{L_pack}. The tests were performed with version {\tt 1.26.06} of \emph{Macaulay2}.

We are confident that the methods and functionalities provided herein may be valuable for further applications. Moreover, the structure of the package allows for future developments, specifically by implementing additional key results and techniques. This task is currently under study by the author.

\subsection*{Acknowledgments} The author thanks Prof. Giancarlo Rinaldo for sharing his original source code, which was later integrated into this package and for his helpful suggestions.
The author also thanks Prof. Marilena Crupi her helpful suggestions.
The author acknowledge support of the GNSAGA group of INdAM (Italy).

\end{document}